\documentclass[11pt]{article}

\usepackage{latexsym}
\usepackage{amsfonts}
\usepackage{enumerate}
\usepackage{multicol}
\usepackage{graphicx}
\usepackage{amssymb}
\usepackage{amsmath}
\usepackage{epic}

\usepackage{amsthm}

\begin{document}
\parindent=0.2in
\parskip 0in

\begin{flushright}

{\huge {\bf Report on the 61st Annual International Mathematical Olympiad}}

\vspace*{.2in}

{\Large B\'ELA BAJNOK} \\

{\small Gettysburg College \\  Gettysburg, PA 17325 \\ bbajnok@gettysburg.edu} 

\vspace*{.2in}

{\Large EVAN CHEN} \\

{\small Massachusetts Institute of Technology \\  Cambridge, MA 02139 \\ evan@evanchen.cc} 

\end{flushright}

\vspace*{.2in}


The International Mathematical Olympiad (IMO) is the world's leading mathematics competition for high school students, and is organized annually by different host countries.  The competition consists of three problems each on two consecutive days, with an allowed time of four and a half hours both days.  In recent years, more than one hundred countries send teams of up to six students to compete.  

The 61st IMO was to take place in July of 2020 in St.~Petersburg, Russia.  Due to the COVID-19 pandemic, however, the competition was postponed and moved remotely where students could take the exam in their home countries.  The competition was ultimately given on September 21 and 22, 2020, at universally coordinated times (which meant that students in the Americas had to start their work in the middle of the night).    

The members of the US team are chosen during the Math Olympiad Program (MOP) each year, a year-long endeavor organized by the MAA's American Mathematics Competitions (AMC) program.  Students gain admittance to MOP based on their performance on a series of examinations, culminating in the USA Mathematical Olympiad (USAMO).  A report on the 2020 USAMO can be found in the April 2021 issue of this {\em Magazine}; more information on the American Mathematics Competitions program can be found on the site https://www.maa.org/math-competitions.  

The members of the 2020 US team were  Quanlin Chen (11th grade, Princeton International School of Mathematics and Science, Princeton, NJ); Gopal Krishna Goel (11th grade, Krishna Homeschool, Portland, OR); Tianze Jiang (11th grade, Princeton International School of Mathematics and Science, Princeton, NJ); Jeffrey Kwan (12th grade, Harker Upper School, San Jose, CA); Luke Robitaille (10th grade, Robitaille Homeschool, Euless, Texas); and William Wang (West Windsor-Plainsboro High School North, Plainsboro Township, NJ).  Chen, Robitaille, and Wang each earned Gold Medals, and Goel, Jiang, and Kwan each earned Silver Medals.  In the unofficial ranking of countries, the United States finished third after China (first) and Russia (second).

Below we present the problems and solutions of the 61st IMO.  The solutions we feature here are those of the indicated students, as edited by the two authors of this report.

\bigskip

\noindent {\bf Problem 1} {\em proposed by Dominik Burek, Poland.}  Consider the convex quadrilateral $ABCD$.  The point $P$ is in the interior of $ABCD$.  The following ratio equalities hold:
$$ \angle PAD:\angle PBA:\angle DPA
	= 1:2:3
	= \angle CBP:\angle BAP:\angle BPC.$$
Prove that the following three lines meet in a point: the internal bisectors of angles $\angle ADP$ and $\angle PCB$ and the perpendicular bisector of segment $AB$.

\bigskip

\noindent {\em Solution by Isabella Quan, 10th grade, Westlake High School, Austin, TX.}  Let $O$ be the circumcenter of $\triangle PAB$.  We will prove that each of the three lines specified in the problem go through $O$.  Since $AO=BO$, this clearly holds for the perpendicular bisector of $AB$.  Below we prove that $O$ lies on the internal bisector of $\angle PCB$; the claim about the internal bisector of $\angle ADP$ can be verified similarly.

\begin{center}
	\includegraphics{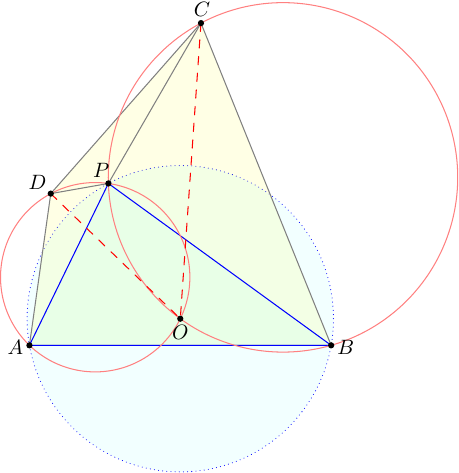}
\end{center}

Note that $$\angle PCB = \pi  - (\angle CBP + \angle BPC)$$ and, by the given condition, $$\angle CBP + \angle BPC = 2\angle BAP.$$  But by the inscribed angle theorem,  $\angle BOP = 2\angle BAP$ as well, and thus $\angle PCB$ and  $\angle BOP$ are supplementary, which implies that quadrilateral $BOPC$ is cyclic.  Because $OP = OB$, the inscribed angles $\angle PCO $ and $ \angle OCB$ in the circumcircle of quadrilateral $BOPC$ are equal, and hence $O$ lies on the angle bisector of $\angle PCB$, as claimed.

\bigskip

\noindent {\bf Problem 2} {\em proposed by Stijn Cambie, Netherlands}.  The real numbers $a, b, c, d$
are such that $a\geq b\geq c\geq d>0$ and $a+b+c+d=1$.
Prove that
$$(a+2b+3c+4d) a^a b^b c^c d^d < 1.$$

\bigskip

\noindent {\em Solution by Luke Robitaille, 10th grade, Robitaille Homeschool, Euless, TX.}  We will prove that, for real numbers $a$, $b$, $c$, and $d$ satisfying the conditions of the problem, we have:
\begin{eqnarray} \label{1}
a^ab^bc^cd^d & \leq & a^2+b^2+c^2+d^2
\end{eqnarray} 
and
\begin{eqnarray}\label{2}
(a+2b+3c+4d)(a^2+b^2+c^2+d^2) & < & 1;
\end{eqnarray} 
these two inequalities together clearly  imply the claim of the problem.

Our first inequality follows from the Weighted AM-GM Inequality, as we now explain.  Let us recall that for a positive integer $n$ and nonnegative real numbers $x_1, \dots ,x_n,$ $w_1, \dots, w_n$ with $w=w_1+ \cdots + w_n >0$, the {\em weighted arithmetic mean} of $x_1, \dots ,x_n$ with corresponding weights $w_1, \dots, w_n$ is defined as $(w_1x_1 + \cdots + w_nx_n)/w,$ and their {\em weighted geometric mean} is defined as
$\sqrt[w]{{x_1}^{w_1} \cdot \cdots \cdot {x_n}^{w_n}}.$   The {\em Weighted AM-GM Inequality} then states that the weighted arithmetic mean is always at least as much as the weighted geometric mean; that is: 
\begin{eqnarray*}
\sqrt[w]{{x_1}^{w_1} \dots {x_n}^{w_n}} & \leq & (w_1x_1 + \cdot \cdots \cdot + w_nx_n)/w.
\end{eqnarray*}
(The case when all weights equal 1 yields the regular AM-GM Inequality.)  Applying the Weighted AM-GM Inequality for $n=4$, $x_1=w_1=a$, $x_2=w_2=b$, $x_3=w_3=c$, and $x_4=w_4=d$, and noting that $a+b+c+d=1$, readily gives (\ref{1}).  

To establish (\ref{2}), we prove more generally that for real numbers $a$, $b$, $c$, and $d$ with $a \geq b \geq c \geq d>0$, we have 
\begin{eqnarray}\label{3}
(a+2b+3c+4d)(a^2+b^2+c^2+d^2) & < & (a+b+c+d)^3.
\end{eqnarray} 

In order to do so, we transform to new variables $r,s,t,u$, defined as $r=a-b$, $s=b-c$, $t=c-d$, and $u=d$; our assumptions imply that $r,s,t \geq 0$ and $u>0$. 
Since $a=r+s+t+u$, $b=s+t+u$, $c=t+u$, and $d=u$, we have
\begin{eqnarray*}
a+2b+3c+4d & = & r+3s+6t+10u
\end{eqnarray*} 
and 
\begin{eqnarray*}
a^2+b^2+c^2+d^2 & = & (r+s+t+u)^2+(s+t+u)^2+(t+u)^2+u^2 \\
& = & r^2+2s^2+3t^2+4u^2+2rs+2rt+2ru+4st+4su+6tu,
\end{eqnarray*}  and thus we get
\begin{eqnarray*}
(a+2b+3c+4d)(a^2+b^2+c^2+d^2) & = & r^3 + 6 s^3  + 18 t^3  + 40 u^3 +    \\
& &  5 r^2 s + 8 r^2 t + 12 r^2 u + 8 r s^2 + 24 s^2 t + 32 s^2 u +          \\
& &  15 r t^2 + 33 s t^2 + 66 t^2 u + 24 r u^2 + 52 s u^2 + 84 t u^2 +         \\
& & 22 r s t  + 30 r s u + 38 r t u + 82 s t u.   
\end{eqnarray*} 
Furthermore,
\begin{eqnarray*}
(a+b+c+d)^3 & = &(r+2s+3t+4u)^3 \\
& = & r^3 + 8 s^3 + 27 t^3   + 64 u^3 + \\
&  & 6 r^2 s + 9 r^2 t + 12 r^2 u  + 12 r s^2  + 36 s^2 t + 48 s^2 u + \\
&  & 27 r t^2 + 54 s t^2   + 108 t^2 u + 48 r u^2 + 96 s u^2 + 144 t u^2 + \\
&  & 36 r s t + 48 r s u + 72 r t u  + 144 s t u .
\end{eqnarray*} 

Making term-wise comparisons, we can verify that 
$$(a+2b+3c+4d)(a^2+b^2+c^2+d^2) \leq (a+b+c+d)^3,$$ and since $u>0$ implies that $40 u^3 < 64 u^3$, we see that, in fact, strict inequality holds, which proves (\ref{3}).  This completes our proof.

\bigskip

\noindent {\bf Problem 3} {\em proposed by Milan Haiman (Hungary) and Carl Schildkraut (USA)}.  There are $4n$ pebbles of weights $1, 2, 3, \dots, 4n$.  Each pebble is colored in one of $n$ colors, and there are four pebbles of each color.  Show that we can arrange the pebbles into two piles, so that the total weights of both piles are the same and each pile contains two pebbles of each color.

\bigskip

\noindent {\em Solution by Alex Zhao, 8th grade, Kamiakin Middle School, Kirkland, WA.}  
We begin by pairing off the pebbles so that each pair sums to the same value;
namely, pairing $1$ and $4n$, $2$ and $4n-1$, $3$ and $4n-2$,
\dots, so that each pair adds up to $4n+1$.
Then, it is sufficient to split these $2n$ \emph{pairs} into two groups of $n$ pairs each,
such that each pile contains two pebbles of each color.

We separate the pebbles by color, and put pebbles of the same color into a box of their own,
and then draw a line between the two pebbles of each pair.
(Some of the lines could start and end in the same box.)
The problem can then be rephrased as follows:
we wish to color each of the lines either blue or green,
such that each box has two pebbles which are endpoints of blue lines,
and two pebbles which are endpoints of green lines.
(These new colors are not related to the colors in the original problem.
They correspond to whether we put the pair in one pile or the other.)
An example with $n=5$ is illustrated below.
\begin{center}
	\includegraphics{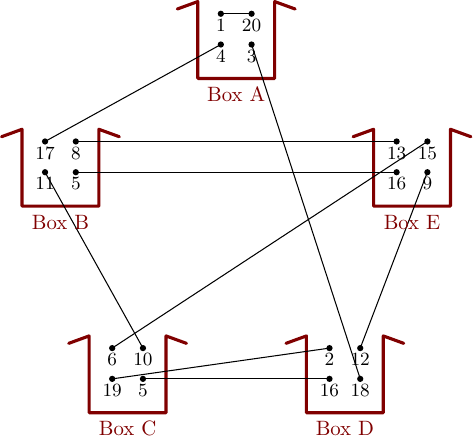}
\end{center}

We can capture this information in a single multigraph $G$.
(A multigraph is like a graph, but edges may be repeated;
we also allow an edge to join a vertex to itself.)
Each box is a vertex of $G$,
and each line is an edge of the multigraph.
Hence, this graph will have $n$ vertices and degree $4$ at every vertex.
Shown below is the multigraph $G$ corresponding to the example we drew before.
(The circle that passes through $A$ denotes the edge from $A$ to itself.)

\begin{center}
	\includegraphics{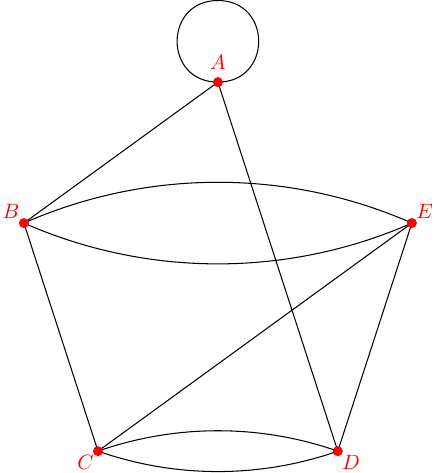}
\end{center}

Consider any connected component of this multigraph.
Since the degree of every vertex is even,
the component has an \emph{Eulerian cycle} ---
a cycle passing through every edge.
(In the above figure, one example of an Eulerian cycle is $AABCDEBECDA$.)

So, consider any connected component and one of its Eulerian cycles.
Since the degree of each vertex is 4,
the sum of all degrees is a multiple of 4,
and so the number of edges (half the sum of the degrees) is even.
Therefore we can alternately color the edges in the Eulerian cycle blue and green.
Shown below is an alternating coloring of $AABCDEBECDA$ in this way.
\begin{center}
	\includegraphics{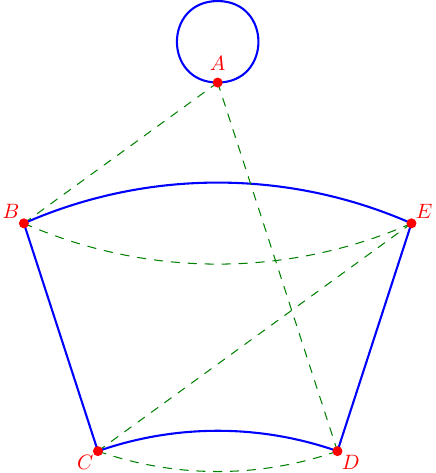}
\end{center}

We can now easily see that this coloring works; indeed:
\begin{itemize}
	\item If the vertex has one self-loop,
		the two outgoing edges must be colored the same,
		and differently from the self-loop.
	\item If the vertex has no self-loops, the cycle visits the vertex twice,
		and the two edges going in and out of each visit
		are colored differently from each other.
	\item If the vertex has two self-loops (which can only happen
		if the connected component has only a single vertex)
		then those two self-loops are colored differently.
\end{itemize}
Having found the desired coloring, the problem is solved.
The figure below shows the translation
of our coloring example into a valid partition of the pairs.
\begin{center}
	\includegraphics{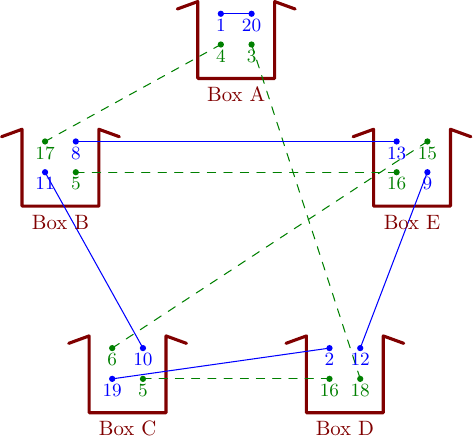}
\end{center}

\bigskip

\noindent {\bf Problem 4} {\em proposed by Tejaswi Navrinarekallu, India}.  Let $n > 1$ be an integer.  There are $n^2$ stations on a slope of a mountain, all at different altitudes.  Each of two cable car companies, $A$ and $B$, operates $k$ cable cars; each cable car provides a transfer from one of the stations to a higher one (with no intermediate stops).  The $k$ cable cars of $A$ have $k$ different starting points and $k$ different finishing points, and a cable car that starts higher also finishes higher.  The same conditions hold for $B$.  We say that two stations are linked by a company if one can start from the lower station and reach the higher one by using one or more cars of that company (no other movements between stations are allowed).
Determine the smallest positive integer $k$ for which one can guarantee that there are two stations that are linked by both companies.

\bigskip

\noindent {\em Solution by Ryan Li, 10th grade, Solon High School, Solon, OH}.
We claim that the answer is $k=n^2-n+1$.

We label the stations by $0$, $1$, $2$, \dots, $n^2-1$
in ascending order of altitude.
Let $G_A$ be a simple graph with $n^2$ vertices defined as follows:
\begin{itemize}
	\item The vertices are the $n^2$ stations.
	\item Two vertices are joined by an edge
		if and only if they are joined by a cable car from company $A$.
\end{itemize}
The graph $G_B$ is defined in an analogous way,
with the connections from company $B$ instead of company $A$.

\begin{quote}
	\textbf{Lemma.} 
	Two stations are linked by $A$ (as defined in the problem statement)
	if and only if they are connected in the graph $G_A$.
\end{quote}

\begin{proof}
	Obviously, being linked implies being connected.
	It remains be verified that the converse is true.

	Suppose we have a path from station $i$ to station $j$ in $G_A$,
	where $i < j$, given by labels $i = a_0 \to a_1 \to  \cdots \to a_{m-1} \to a_m = j$.
	We claim that \[ i < a_1 < a_2 < \cdots < a_{m-1} < j. \]

	Indeed, note that for $r=1,2,\dots,m-1$,
	each $a_r$ is connected to $a_{r-1}$ and $a_{r+1}$ by cable cars .
	Since the starting points and ending points of each cable car are unique,
	$a_{r-1}$ and $a_{r+1}$ cannot be both higher or lower than $a_r$,
	so either $a_{r-1} < a_r < a_{r+1}$ or $a_{r-1} > a_r > a_{r+1}$ for each $r$.

	Since $i < j$, this implies $i = a_0 < a_1$, $a_1 < a_2$, \dots,
	$a_{m-1} < a_m = j$, as claimed.
\end{proof}

Now we prove that $k=n^2-n+1$.
First, we show that there exists a construction for $k=n^2-n$
such that no two stations are linked by both company $A$ and company $B$.
Consider the situation in which for every integer $0 \leq i \leq n-1$,
company $A$ links together $i, i+n, i+2n, \dots, i+n^2-n$ in a series
and company $B$ links together stations $ni, ni+1, ni+2, \dots, ni+n-1$.
This gives a total of $n^2-n$ cable cars for each company.
Depicted below is the case for $n=4$,
with cars for $A$ in bold, and cars for $B$ dashed.
\begin{center}
	\includegraphics{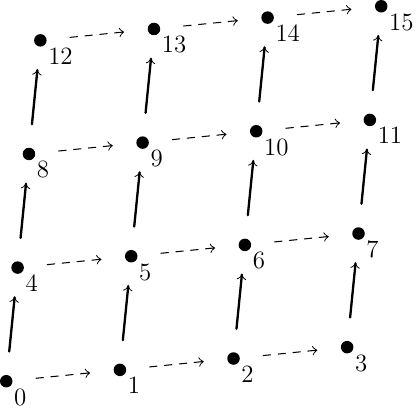}
\end{center}
It is clear that no car is linked by both companies:
since any two cars linked by $A$ have indices differing by at least $n$,
while any two cars linked by $B$ have indices differing by at most $n-1$.


Finally, we prove that $k = n^2 - n + 1$ suffices.
In the earlier proof of the lemma,
we saw that a path from a vertex $i$ to a vertex $j$ must in fact
pass through strictly increasing indices.
This implies that $G_A$ and $G_B$ are acyclic --- i.e., they are forests.

Now since $G_A$ has $n^2$ vertices and $k$ edges but no cycles,
the number of connected components in $G_A$ is exactly
\[ n^2 - k = n^2 - (n^2 - n+ 1) =  n - 1. \]
So by the pigeonhole principle, some connected component $S$ of $G_A$ has
at least $\left\lceil \frac{n^2}{n-1} \right\rceil = n + 1$ vertices.
On the other hand, $G_B$ also has exactly $n - 1$ connected components.
So by the pigeonhole principle again, there are two vertices in $S$ which lie
in the same connected component of $G_B$.
These two stations are therefore linked by both companies.

\bigskip

\noindent {\bf Problem 5} {\em proposed by Oleg Ko\v{s}ik, Estonia.}  A deck of $n > 1$ cards is given.  A positive integer is written on each card.
The deck has the property that the arithmetic mean of the numbers on each pair of cards is also the geometric mean of the numbers on some collection of one or more cards.  For which $n$ does it follow that the numbers on the cards are all equal?

\bigskip

\noindent {\em Solution by Pravalika Putalapattu, 10th grade, Thomas Jefferson High School for Science and Technology, Alexandria, VA.}  We will prove that this assertion is true for all $n$.  

Proceeding indirectly, let us assume that there exists some positive integer $n$ and a set of cards labeled $a_1, a_2,\ldots, a_n$ that satisfy the above property but are not all equal.  We may assume without loss of generality that our cards are in decreasing order of their labels:
\begin{equation*}
    a_1\geq a_2 \geq \cdots \geq a_n.
\end{equation*}
Note that multiplying all card values by some constant will multiply all arithmetic and geometric means by the same value.  Thus, the property is preserved if we replace each label $a_i$ with $a_i/d$, where $d=\gcd(a_1, a_2,\ldots, a_n)$.  Therefore, we may assume without loss of generality that the $n$ labels are relatively prime.

If $a_1=1$, then all $a_i$ must equal 1, violating our assumption that not all $a_i$ are equal, thus, $a_1\geq 2$.  Let $p$ be a prime divisor of $a_1$. Since $\gcd(a_1, a_2, \ldots, a_n)=1$, there must be some $a_i$ not divisible by $p$.  Let $k$ be the smallest index such that $p$ does not divide $a_k$. Note that $a_k < a_1$.

Now, consider the arithmetic mean of $a_1$ and $a_k$, which we know is greater than $a_k$. This must equal some geometric mean, and thus 
\begin{equation}
    \frac{a_1+a_k}{2}=\sqrt[m]{a_{i_1}  \cdots  a_{i_m}} \nonumber
\end{equation}
for some positive integer $m$ and indices $i_1, \dots, i_m$.  In particular, $$\sqrt[m]{a_{i_1}  \cdots a_{i_m}} > a_k.$$  

Next, observe that, since $a_k$ is not divisible by $p$, none of $a_{i_1}, \dots, a_{i_m}$ can be divisible by $p$ either.   Since $k$ is the smallest index such that $p$ does not divide $a_k$, each of $i_1, \ldots, i_m$ must be at least $k$, and therefore each of  $a_{i_1}, \ldots, a_{i_m}$ must be at most $a_k$.  But then their geometric mean is also at most $a_k$, resulting in a contradiction, as desired.  

\bigskip

\noindent {\bf Problem 6} {\em proposed by Ting-Feng Lin and Hung-Hsun Hans Yu, Taiwan}.  Consider an integer $n > 1$ and a set $\mathcal S$ of $n$ points
in the plane such that the distance between any two different points in $\mathcal S$ is at least $1$.  Prove that there is a line $\ell$ separating $\mathcal S$ such that the distance from any point of $\mathcal S$ to $\ell$
is at least $cn^{-1/3}$ for some positive constant $c$.
(A line $\ell$ separates a set of points $S$ if some segment joining two points in $\mathcal S$ crosses $\ell$.)

\bigskip

\noindent {\em Solution by Gopal Krishna Goel, 11th grade, Krishna Homeschool, Portland, OR, and Jaedon Whyte, 10th grade, Archimedean Upper Conservatory, Miami, FL.}  
We will show that there exists a line $\ell$ separating $\mathcal{S}$ such that the distance from any point in $\mathcal{S}$ to $\ell$ is at least $0.01\cdot n^{-1/3}$.

Let $[\mathcal{P}]$ denote the area of a region $\mathcal{P}$ in the plane. We will need the following lemma.

\begin{quote}
\textbf{Lemma}.
Let $\mathcal{R}$ be an $a\times b$ rectangle in the plane, where $a,b\ge\tfrac{1}{2}$. Then $|\mathcal{S}\cap\mathcal{R}|\le 20[\mathcal{R}]$.
\end{quote}
\begin{proof}
For each point $P$, let $\mathcal{D}_P$ denote the open disk of radius $\tfrac{1}{2}$ centered at $P$. Since $PQ \ge 1$ for all distinct $P,Q\in\mathcal{S}$, their respective disks $\mathcal{D}_P$ and $\mathcal{D}_Q$ cannot intersect.

Now, for any $P\in\mathcal{S}\cap\mathcal{R}$, we see that $\mathcal{D}_P\subseteq\mathcal{R}'$, where $\mathcal{R}'$ is an $(a+1)\times(b+1)$ rectangle with each side of $\mathcal{R}'$ a distance exactly $1/2$ from a side of $\mathcal{R}$.
For the total area covered by the disks we have:
\[
	\sum_{P\in\mathcal{S}\cap\mathcal{R}}[\mathcal{D}_P]\le[\mathcal{R}']
	= |\mathcal{S}\cap\mathcal{R}|\cdot\frac{\pi}{4}
	= (a+1)(b+1).\]
Since $a,b\ge\tfrac{1}{2}$, we see that $\frac{(a+1)(b+1)}{ab} = \left(1+\frac{1}{a}\right)\left(1+\frac{1}{b}\right) \leq 3\cdot 3 = 9$, so:
\[|\mathcal{S}\cap\mathcal{R}|\le \frac{4}{\pi}(a+1)(b+1)\le\frac{36}{\pi}ab<20[\mathcal{R}],\]
as desired.
\end{proof}
Let $A$ and $B$ denote the two points in $\mathcal{S}$ such that $AB$ is maximized (if there are multiple pairs, pick one arbitrarily), and let $r=AB$ be this maximum value. Let $\mathcal{T}$ denote the set of points $P$ in the plane such that $\max(AP,BP) \le r$. We see that $\mathcal{T}$ is simply the intersection of the disk centered at $A$ with radius $r$ and the disk centered at $B$ with radius $r$, as shown below. By definition, we must have $\mathcal{S}\subseteq\mathcal{T}$.

\begin{center}
	\includegraphics{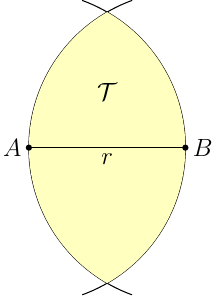}
\end{center}

We will first resolve the case $r\ge n^{2/3}$. Project all the points of $\mathcal{S}$ onto $AB$, and suppose they are $A=P_1,P_2,P_3,\ldots,P_n=B$, appearing in that order (it is possible that we may have $P_i=P_{i+1}$ for some values of $i$). Now, since
\[P_1P_2+P_2P_3+\cdots+P_{n-1}P_n = r,\]
we see that there exists some $i$ such that
\[P_iP_{i+1}\ge \frac{r}{n-1}>r/n\ge n^{-1/3}.\]
Take $\ell$ to be the perpendicular bisector of $P_iP_{i+1}$. We see that the distance from any point $P\in\mathcal{S}$ to $\ell$ is simply the distance from the projection of $P$ onto $AB$ to the midpoint of $P_iP_{i+1}$, which is at least 
\[\frac{1}{2}P_iP_{i+1}>\frac{1}{2}n^{-1/3},\]
as desired. This resolves the case $r\ge n^{2/3}$.

We now turn our attention to the case $r<n^{2/3}$. Let $XY$ be a chord of the circle centered at $B$ with radius $r$ such that $XY\perp AB$ and the distance from $A$ to $XY$ is $\tfrac{1}{2}$. Let $M$ be its midpoint. Since $r=AB\ge 1$, we see that $X$ and $Y$ are on the arc of the circle centered at $B$ with radius $r$ that is a part of $\mathcal{T}$.
The line $XY$ divides $\mathcal T$ into two parts;
let $\mathcal{T}'$ be the part containing $A$.
\begin{center}
	\includegraphics{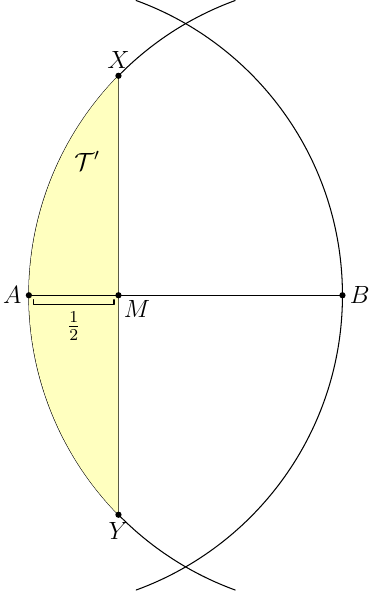}
\end{center}
Note that $\mathcal{T}'$ can be contained in a rectangle of size $\tfrac{1}{2}\times XY$. By the Pythagorean theorem on $\triangle{XMB}$, we have
\[XY = 2\sqrt{r^2-\left(r-\frac{1}{2}\right)^2} = 2\sqrt{r-\frac{1}{4}}<2\sqrt{r}<2n^{1/3},\]
so by the lemma above, we see that
\[|\mathcal{T}'\cap\mathcal{S}|\le 20\cdot\frac{1}{2}\cdot|XY| < 20 n^{1/3}.\]
Similar to before, let $Q_1=A,Q_2,Q_3,\ldots,Q_{|\mathcal{T}'\cap\mathcal{S}|}$ denote the projections of the points of $\mathcal{T}'\cap\mathcal{S}$ onto segment $AB$, in that order.

If $A$ is the only point of $S$ in $\mathcal T'$,
then we take our line as the perpendicular bisector of line $AM$.
So assume $|\mathcal{T'} \cap \mathcal{S}| \ge 2$.
Since all points $Q_i$ lie in a segment of length $\tfrac{1}{2}$, so there must be some $i$ such that
\[Q_iQ_{i+1}\ge\frac{\tfrac{1}{2}}{|\mathcal{T}'\cap\mathcal{S}|-1}>\frac{1}{2|\mathcal{T}'\cap\mathcal{S}|}>\frac{1}{40}n^{-1/3}.\]
We may now take $\ell$ to be the perpendicular bisector of $Q_iQ_{i+1}$, and as before, the distance from any point of $\mathcal{S}$ to $\ell$ is at least $\tfrac{1}{80}n^{-1/3}>0.01 n^{-1/3}$, as desired. This completes the solution.

\vspace*{.3in}

\noindent {\bf B\'ELA BAJNOK} (MR Author ID: 314851) is a Professor of Mathematics at Gettysburg College and the Director of the American Mathematics Competitions program of the MAA.

\noindent {\bf EVAN CHEN} (MR Author ID: 1158569) is a graduate student at the Massachusetts Institute of Technology and the Assistant Academic Director of the Math Olympiad Program of the MAA.

\end{document}